\documentclass[12pt,leqno]{amsart}
\usepackage{amssymb}
\usepackage{amsmath,amssymb}
\usepackage{graphicx} %% newly added
\theoremstyle{plain}

\theoremstyle{definition}

\numberwithin{equation}{section}

\newcommand{\R}{\mathbb{R}}

\newcommand{\www}{\widetilde}

\newcommand{\ppp}{\partial}
\newcommand{\ooo}{\overline}

\title
[]
{Conditional stability for an inverse source problem and 
an application to the estimation of air dose rate radioactive 
substances by drone data}

\begin{document}

\author[Y. Chen, J. Cheng, G.Floridia, Y. Wada, M.Yamamoto]
{Yu Chen $^1$, Jin Cheng $^2$,  Giuseppe Floridia $^3$, 
Youichiro Wada $^4$, Masahiro Yamamoto $^5$} 
\thanks{$^1$ Yu Chen, $^2$ Jin Cheng, School of Mathematical Sciences, 
200433 Shanghai China,\\
email: {\tt yu\_{}chen11@fudan.edu.cn $^1$, jcheng@fudan.edu.cn $^2$}, \\ 
$^3$ Giuseppe Floridia, Department of Mathe\-matics and Applications \lq\lq 
R. Caccioppoli, University of Naples Federico II,
80126 Naples, Italy, email: {\tt giuseppe.floridia@unina.it \& 
floridia.giuseppe@icloud.com},\\
$^4$ Youichiro Wada,
Isotope Science Center, The University of Tokyo, Tokyo, 113-0032, Japan,
email: {\tt wada-y@lsbm.org}, \\
$^5$ Masahiro Yamamoto, Department of Mathematical Sciences, The University
of Tokyo, Komaba, Meguro, Tokyo 153, Japan,
e-mail: {\tt myama@ms.u-tokyo.ac.jp}.}

\baselineskip 18pt

\begin{abstract}  
\medskip
\noindent
%{\bf  Keywords:} Inverse source problems, 
%fractional diffusion equation, reconstruction, well-posedness, stability 
%estima%te.\\
%\medskip
%\noindent
%{\bf Mathematics subject classification 2010 :} 35R30, 	35R11.
%
We consider the density field $f(x)$ generated by a volume source $\mu(y)$ in
$D$ which is a domain in $\R^3$.  For two disjoint segments $\gamma, \Gamma_1$
on a straight line in $\R^3 \setminus \ooo{D}$, 
we establish a conditional stability estimate of H\"older type in determining
$f$ on $\Gamma_1$ by data $f$ on $\gamma$.
This is a theoretical background for real-use solutions for the 
determination of air dose rates of radioactive
substance at the human height level by high-altitude data. 
The proof of the stability estimate is based on the harmonic extension and 
the stability for line unique continuation of a harmonic function.
\end{abstract}
\maketitle
\section{Motivation}

The Fukushima Daiichi Nuclear Disaster in March 2011 has released 
radioactive substances such as cesium-137 into environments.  
In particular, radioactive substances fall on the ground and some of them
diffuses in the soil.  These radioactive substances are a source and generate 
the air dose rate of radiation, which can be considerd as   
influences to inhabitants.
Possible sources for such air dose rates and possible sources are
in the air and on the ground, in the underground.
The sources in the air could be created by floating radioactive substances, but
we can assume that such sources in the air can be neglected after 9 years 
passed after the disaster in 2011 and the susbstances have sufficiently 
diffused in the air to a very low level.
Thus we assume that the air dose rates can be generated by sources on the 
ground and the underground.  As for related works on the several kinds of 
sources, we refer to Malins, Okumura, Machida, Takeyama and Saito
\cite{MOMTS}, Saito \cite{S} for example.
 
Moreover these substances on and in the ground can move for example 
by refloating on dry days and 
may run-off into rivers, so that the sources may not be stationary.
However for a feasible model, we can assume that we can neglect also the time
change of sources.  In some areas of Fukushima prefecture, inhabitants
have already returned to daily lives, and the estimation of the air dose
rate of radioactive substances in towns and villages is crucial for the sake 
of the health of them.  The air dose rate at the human height level (e.g., 1m) 
is considered as very influential factor to the health through exposures by 
breathing the air.
The observation of air dose rates can be done by drones 
containing measurement equipments.
However direct observations of air dose rates at the human height level is 
not realistic because of many obstacles against the flights of drones such as
houses and other many artificial structures such as walls.
Moreover we do not a priori know the source, that is,
the density disribution of radioactive substances on the ground and 
in the shallow underground.
Thus, not knowing sources themselves, we can observe only data 
at higher altitude (e.g., 30m $\sim$ 50m).

Therefore our task is an inverse problem where
we are required to determine the dose rate at the human height leval by means 
of high-altitude data.
This motivates our theoretical research for the correponding inverse problem
and our theoretical result shows not bad stability under adequete conditions
and suggests how to control flight orbits for a reasonable accuracy in 
solving an inverse problem.

\section{Mathematical model and the main result.}

We describe a mathematical model, by which we can discuss our inverse problem 
from a general point of view.

Let $D \subset \R^3$ be an open domain.
Henceforth we set $x = (x_1,x_2,x_3)$, $y= (y_1,y_2, y_3)$ and 
$$
r_{xy} = ((x_1-y_1)^2 + (x_2-y_2)^2 + (x_3-y_3)^2)^{\frac{1}{2}}.
$$
We set 
$$
f(x) := \int_D \frac{\mu(y)}{r_{xy}^2} dy, \quad x \in \R^3\setminus 
\ooo{D}.
                                                            \eqno{(2.1)}
$$
Here $4\pi\mu(y)$, $y \in D$ describes the density of e.g., 
radioactive substances.  We consider $f(x)$ as the air dose rate at 
$x\in \R^3$.
For an arbirarily fixed constant $M>0$, we set 
$$
\mathcal{M} := \{ \mu \in L^{\infty}(D);\, 
\Vert \mu\Vert_{L^{\infty}(D)} \le M\}.
$$
Let $\gamma, \Gamma_1 \subset \R^3$ be sets and $\ooo{\gamma} \cap 
\ooo{\Gamma_1} = \emptyset$.
\\

Our inverse problem is: determine $f\vert_{\Gamma_1}$ by 
$f\vert_{\gamma}$.
\\

Here we note that we assume $\mu\in \mathcal{M}$ but we do not know
$\mu$.

{\bf Example:}\\
Taking into consideration the observation by drones, we  
choose curves $\Gamma_1, \gamma$ on the same curve $\Gamma$ 
which is an orbit of a drone.
For example, we choose 
two segments $\gamma$ and $\Gamma_1$.  
Let $x_3 = 0$ and $x_3 > 0$ correspond to
the flat ground and the underground respectively.  
Let $\www{D} \subset \R^3$ be an 
open domain such that $\www{D} \cap \{ x_3 < 0\} \ne \emptyset$.
Then we set $D = \www{D} \cap \{ x_3 < 0\}$.  In particular, the case of
$\www{D} \subset \{0<x_3<\delta\}$ with small $\delta>0$ and 
$\mu \ne 0$ in $D$, means that the substances concentrate on the ground and 
the shallow underground.
Let $(a,b) \in \R^2$ be a fixed planar location and 
$\delta < h_1 < h_2 < H_1 < H_2$ and $$
\gamma = \{ (a,b,x_3);\, H_1 < x_3 < H_2\}, \quad
\Gamma_1 = \{ (a,b,x_3);\, h_1 < x_3 < h_2\}.    \eqno{(2.2)}
$$
This setting describes that we are requested to determine 
the dose rate at the heights $h_1 \sim h_2$ by high-altitude data
at $H_1\sim H_2$.
\\

Let $L$ be an infinite straight line and let $\Gamma$ be 
a connected component of $L \cap E$.
For an arbitrarily fixed $\delta>0$, we set
$$ 
E = \{ x\in \R^3;\, \mbox{dist}\, (x, D) > \delta\}. 
$$

{\bf Theorem.}
\\
Let $\gamma, \Gamma_1 \subset \Gamma$ be finite segments on 
$\Gamma$ such that 
$\ooo{\gamma} \cap \ooo{\Gamma_1} = \emptyset$.  We assume that 
$$
\mu \in \mathcal{M}.                   \eqno{(2.3)}
$$
Then there exist constants $C>0$ and $\theta \in (0,1)$ such that 
$$
\Vert f\Vert_{L^{\infty}(\Gamma_1)} \le C\Vert f\Vert_{L^{\infty}(\gamma)}
^{\theta}.                         \eqno{(2.4)}
$$
Here $C>0$ depends on $\frac{M}{\delta^2}$, $\gamma, \Gamma_1$, and 
$\theta$ depends only on $\gamma, \Gamma_1$. 
\\

We can expect that $\theta$ is smaller when dist $(\gamma, \Gamma_1)$ is 
larger.
The theorem implies the uniqueness, that is, if $f=0$ on $\gamma$, then 
$f=0$ on $\Gamma_1$.

As can be seen by the proof in Section 3, data $f$ on $\gamma$ cannot 
give any information of $f$ outside of the straight line $L$ where 
$\gamma$ and $\Gamma_1$ are included. 
In other words, we can determine $f$ only in the extended direction of
$\gamma$ by data on $\gamma$, provided that the extended segment is in 
$E$.
In particular, we assume that $0<\delta<h_1<h_2<H_1<H_2$ and recall (2.2).
Then $\gamma, \Gamma_1 \subset \Gamma:= (\delta, \infty)$.  The theorem 
asserts a stability estimate of H\"older type in determining 
$f(x)$, $x\in \Gamma_1$ by $f(x)$, $x \in \gamma$, and the H\"older 
exponent $\theta$ depends only on a geometric configuration of the 
two segments $\gamma$ and $\Gamma_1$.  Our main theorem guarantees a
rather good rate of the stability when we use drone data at high-altitude 
for determining the dose rate at lower-level.

Here we do not discuss the determination of the source $\mu(x)$ itself, and 
we do not know whether we can extract some information of $\mu$ 
from $f\vert_{\gamma}$.  On the other hand, as data we choose 
$f$ in a domain $D_1 \subset \R^3$ such that 
$\ooo{D \cap D_1} = \emptyset$, and Cheng, Pr\"ossdorf and Yamamoto 
\cite{CPY} proves some stability for an inverse problem of 
determining $\mu\vert_{D}$ by $f\vert_{D_1}$ which are not data on 
a segment.  The stability is conditional under assumption that unknown 
$\mu$ satisfies (2.3) and is of logarithmic rate which is much weaker 
than (2.4).
We emphasize that data $f$ on $\gamma$ cannot determine the source
$\mu$ but can determine $f$ on $\Gamma$ with H\"older stability rate which is
much better than the logarithmic rate.

\section{Proof of Theorem}

The proof is based on a harmonic extension of $f$ and 
the line unique continuation for a harmonic function.

{\bf First Step:harmonic extension of $f$}\\
We recall $E = \{x\in \R^3;\, \mbox{dist}\, (x,D) > \delta\}$ with 
arbitrarily fixed $\delta>0$.  We set
$$
G(x,\xi) = \int_D \frac{\mu(y)}{r_{xy}^2+\xi^2} dy, \quad x\in E, \,
\xi \in \R.                                \eqno{(3.1)}
$$
We can directly verify 
$$
\Delta G = \sum_{j=1}^3 \frac{\ppp^2 G}{\ppp x_i^2} 
+ \frac{\ppp^2G}{\ppp \xi^2} = 0 \quad \mbox{in $E\times \R$}.
$$
Clearly 
$$
G(x,0) = f(x), \quad x \in E,           \eqno{(3.2)}
$$
which means that $G$ is a harmonic extension of $f$.

We note that the harmonic extension $G$ defined by (3.1) is applied for the 
determination of
$\mu$ (e.g., Cheng, Pr\"ossdorf and Yamamoto \cite{CPY},
Cheng and Yamamoto \cite{CY2}).

{\bf Second Step:line unique contination}\\

We set 
$$
\www{\Gamma} = \{ (x,0);\, x\in \Gamma\}, \quad
\www{\gamma} = \{ (x,0);\, x\in \gamma\},\quad
\www{\Gamma_1} = \{ (x,0);\, x\in \Gamma_1\} \subset \R^4.
$$
Then $\www{\Gamma}, \www{\gamma}, \www{\Gamma_1} \subset E\times \{0\}$
are lines.  Moreover for any $\mu \in \mathcal{M}$, we can estimate
$$
\vert G(x,\xi)\vert \le \int_D \frac{\vert \mu(y)\vert}
{r_{xy}^2} dy \le \frac{M}{\delta^2}, \quad (x,\xi) \in E \times \R.
                                                 \eqno{(3.3)}
$$
In view of (3.3), we apply the line unique continuation by 
Cheng, Hon and Yamamoto \cite{CHY} (also Cheng and Yamamoto \cite{CY1}), so 
that we can choose constants $C=C\left(\frac{M}{\delta^2},\gamma,\Gamma_1
\right)> 0$ and $\theta=\theta(\gamma,\Gamma_1)\in (0,1)$ such that 
$$
\Vert G\Vert_{L^{\infty}(\www{\Gamma_1})} 
\le C\Vert G\Vert_{L^{\infty}(\www{\gamma})}^{\theta},  \eqno{(3.4)}
$$
that is,
$$
\Vert f\Vert_{L^{\infty}(\Gamma_1)} 
\le C\Vert f\Vert_{L^{\infty}(\gamma)}^{\theta}
$$
by (3.2).  Thus the proof is complete.

\section{Numerical examples}
In the numerical computation, we here propose a method via 
the determination of $\mu(x)$.  More precisely,
we choose $\{p_1,p_2,...,p_M\}\subset \gamma\subset\mathbb{R}^3$ and 
$\{q_1,q_2,...,q_N\}\subset D\subset\mathbb{R}^3$ as collocation points, and 
construct the matrix $\mathbf{A}=[A_{ij}]$ as $A_{ij}=w_j/r^2_{p_iq_j}$ where 
$w_j$ are the volume integral coefficients, i.e.,
\[
\int_{D} \frac{\mu(y)}{|p_i-y|^2}\mathrm{d}V(y)
\sim \sum_{j=1}^N A_{ij}\mu(q_j),\quad i=1,...,M.
\]
Then find the least-norm solution $\tilde{\mu}$ to the equation 
$\mathbf{A}\tilde{\mu}=\mathbf{b}$, where $\mathbf{b}=(f(p_1),...,f(p_M))$. 
Taking into consideration observation data $f_{meas}$ with errors, 
the Tikhonov regularization is 
adopted and $\tilde{\mu}$ is then the minimizer to the cost functional
$$
J(\mu):=\|f(\mu)-f_{meas}\|^2_{L^2(\gamma)}+\alpha\| \mu \|^2_{L^2(D)}.
$$
Here $\alpha > 0$ is a regularizing parameter which we should choose
suitably according to the noise level.
A discretized form is
$$
J(\tilde{\mu})=\sum_{i=1}^M\sum_{j=1}^N(A_{ij}\mu_j-b_i)^2 l_i +\alpha 
\sum_{j=1}^N \mu_j^2  w_j,
$$
where $l_i$ are the line integral coefficients. 
The regularization parameter $\alpha$ is chosen as $\alpha\sim \delta^2$, where $\delta$ is the noise level of observation data measured by the 
$L^2(\gamma)$-norm and we refer to Cheng and Yamamoto \cite{CY3} as for a
choice strategy of $\alpha$.
The approximation solution can be constructed by $\tilde{\mu}$ as 
\[
\tilde{f}(x)=\sum_{j=1}^N \frac{\tilde{\mu}_j w_j}{|x-q_j|^2}.
\]
In the following numerical examples, we assume that 
the source distribution has a compact support in $B:= \{x=(x_1,x_2,x_3)\in 
\R^3; \,\vert x\vert < 1, \, x_3<0\}$.
The exact source distribution is given as 
$\mu(x)=(1+0.5\sin(2\pi \vert x\vert))\exp(-\vert x\vert^2)
\chi_B(x)$ where $\chi_B(x)$ is the characteristic function of $B$. 

\begin{figure}%[H] 
\begin{center}
\includegraphics[width=2.5in]{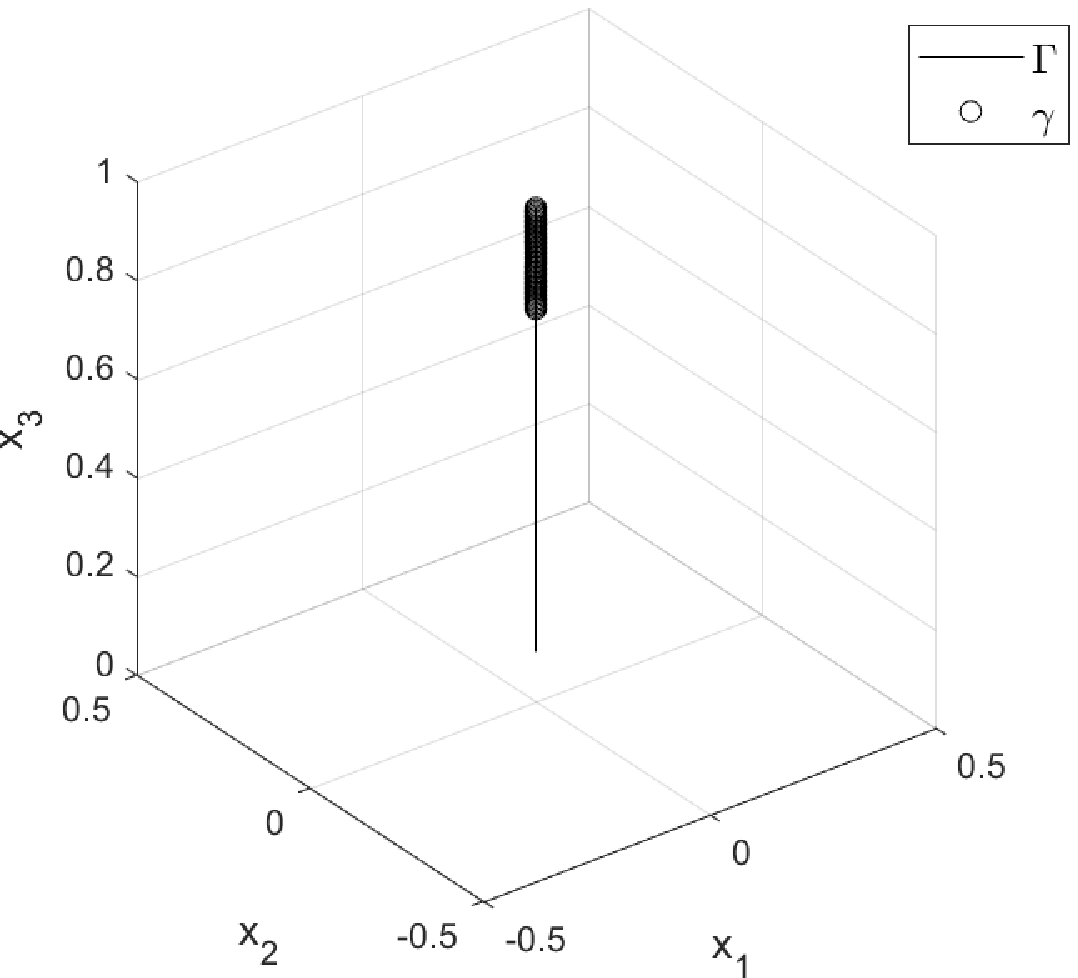}
\includegraphics[width=2.75in]{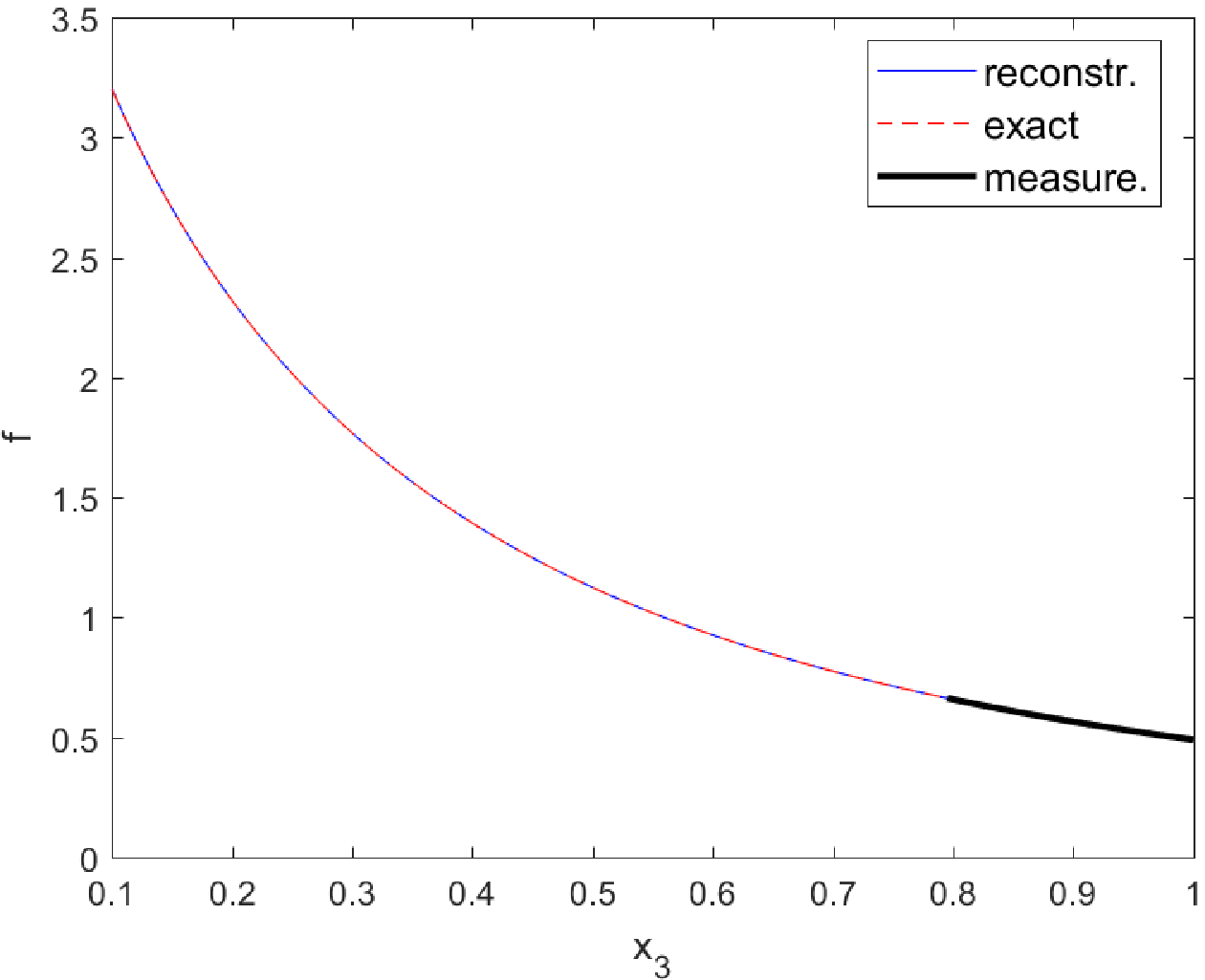}
\end{center}
  \caption{Reconstruction of $f|_{\Gamma_1}$ by $f|_\gamma$ on a line. Left: 
Sketch of the line and measurement points; right: numerical result as a 
function of the height ($x^3$).}
\label{case-1} 
\end{figure}

The present example is to reconstruct the dose rate along the line 
$\Gamma=\{0,0,x_3\}$ with measurement on $\gamma=\{0,0,x_3\},x_3\in [0.8,1]$. 
There are 20 uniformly distributed collocation points on the measured segment 
and $40\times40\times 20$ integral points in $[-1.0,1.0]\times[-1.0,1.0]
\times[-1.0,0]$. Figure \ref{case-1} shows the comparison between 
the reconstructed solution and the exact solution along $\Gamma_1=\{0,0,x_3\},
x_3\in [0.1,1]$ in the case of no observation noises. 
The exact solution are obtained by numerical integration 
and the difference is less than $10^{-4}\|f\|_{C(\Gamma)}$ with finer grid 
resolution.

\begin{figure}%[H] 
\begin{center}
\includegraphics[width=2.75in]{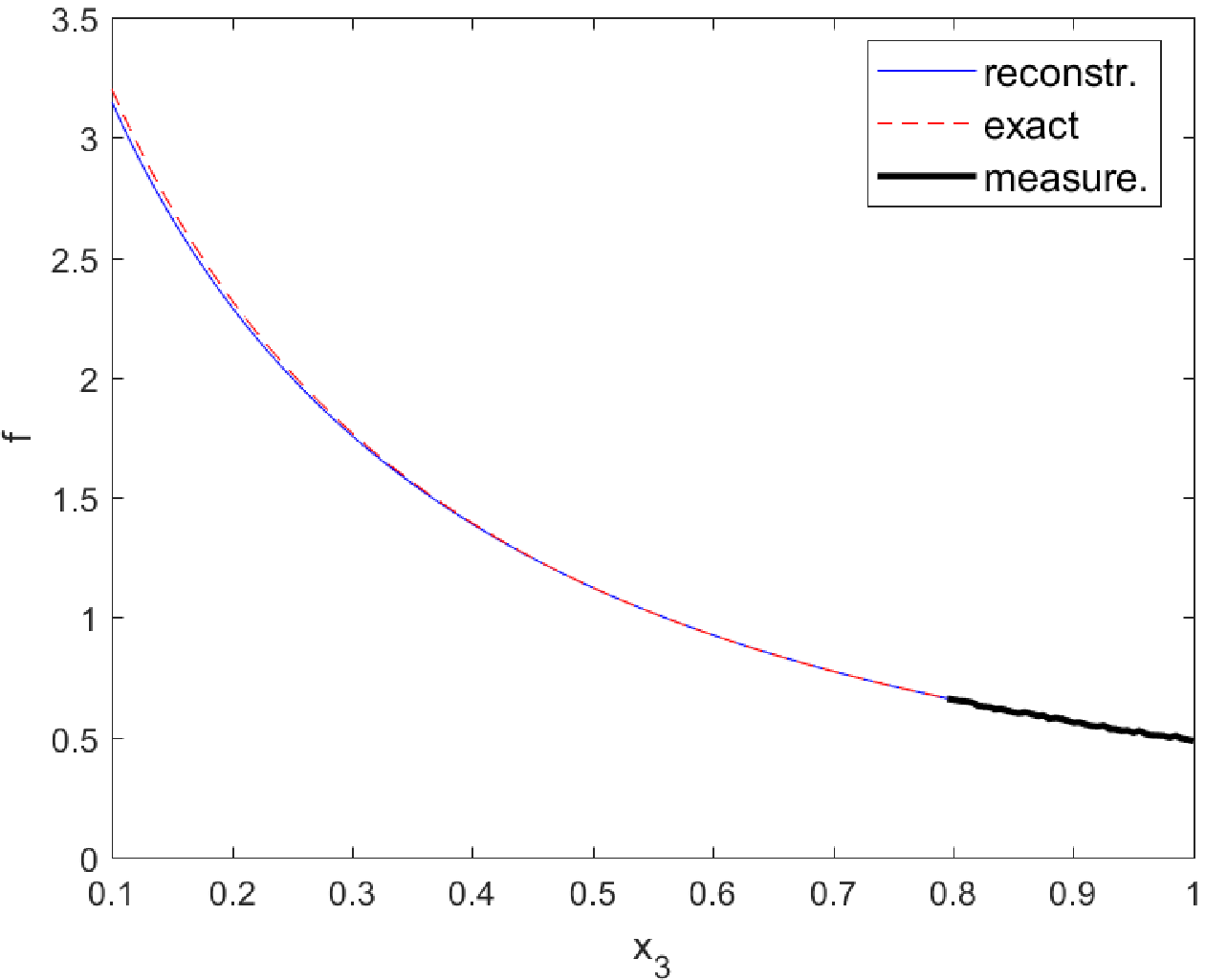}
\includegraphics[width=2.75in]{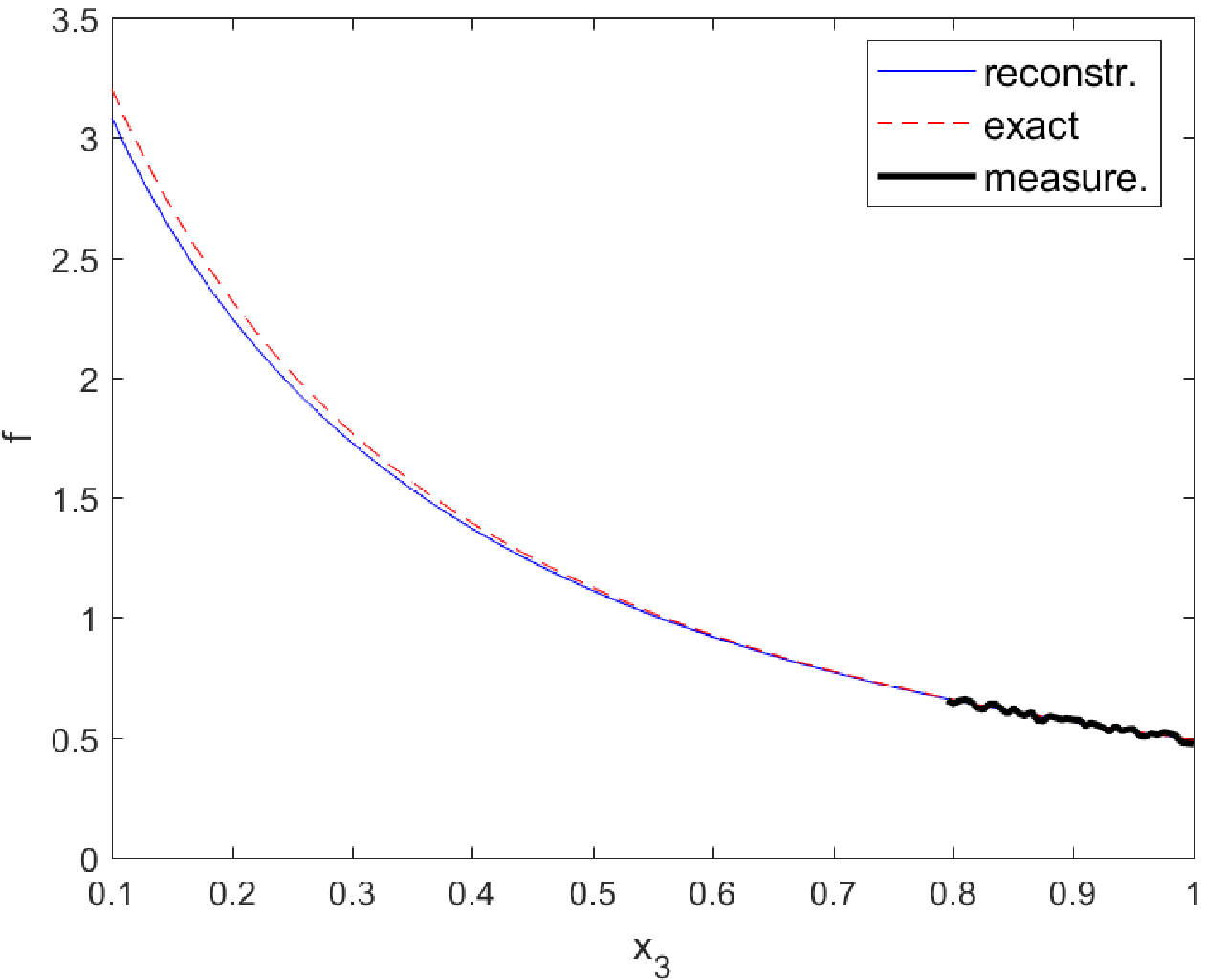}
\end{center}
  \caption{Reconstruction of $f|_{\Gamma_1}$ by $f|_\gamma$ on a line with $1\%$ noize (left) and $3\%$ noize (right) in measurement.}
\label{case-2} 
\end{figure}

In order to further illustrate the performance of the method, a random noise 
is added on the exact observation. %, with the pointwise error within $3\%$ of $\|f\|_{C(\gamma)}$.
Figure \ref{case-2} gives the results with observation error of level $1\%$ 
and $3\%$ respectively, which coincide well with the exact solution with the 
relative error less than $5\%$ at $x_3=0.1$. 

\section{Conclusion}

{\bf 1.}\\
In this paper, we dicuss the estimation of air dose rates of 
radioactive substances at the human height level by high-altitude data
by means by drone.  We establish a H\"older stability estimate for this 
inverse problem on the basis of a mathematical model.
For determination of air dose rates in more directions, 
the main theorem suggests
that a drone orbit should include the corresponding many directions.
More precisely, let $\gamma$ be a curve composed of $N$-segments $\gamma_1, 
..., \gamma_N \subset E$.
Let $\ell_k$, $1\le k \le N$ be the extended straight lines of $\gamma_k$ 
such that $\ell_k$ and $\gamma_k$ are in the same connected component 
of $E$.  The line $\ell_k$ may be a finite segment, an infinite straight line 
or a half straight line.  The theorem means that 
$f\vert_{\gamma}$ determines $f$ on $\cup_{k=1}^N \ell_k$ uniquely.
In other words, when a drone flies zigzag, drone data can 
determine $f$ on all the extended directions of the components of the 
zigzag.  
Moreover the stability is rather good in determining the dose rate 
along the orbit direction.

{\bf 2.}\\
The key of the proof is the harmonic extension and the line unique 
continuation, which asserts stable continuation (3.4) of a harmonic 
function along the segment and is different from usual unique continuation 
of solutions to elliptic equations (e.g., Isakov \cite{I},
Lavrent'ev, Romanov and Shishatskii \cite{LRS}). 
The uniqueness corresponding to (3.4) can be easily proved.  Indeed
let $x(\tau)$ with $\tau \in J$: some open interval, be a parametrization 
of the straight line $\Gamma$.  Since $\Delta G = 0$ in some domain 
$\www{E} \subset \R^4$, the interior real analyticity of a harmonic
function yields that $G$ is real analytic in $(x,\xi) \in \www{E}$.
This analyticity can be proved also by the representation (3.4) of $G$.
Therefore $G(\tau):= G(x(\tau),0)$ is analytic in $\tau \in J$.
Consequently with some non-empty open interval $J_0 \subset J$, the
unicity theorem of an analytic function implies that
if $G(\tau) = 0$ for $\tau \in J_0$, then $G(\tau) = 0$ for $\tau \in J$.

For more details on the line unique continuation, we refer to 
\cite{CHY} and \cite{CY1}

{\bf 3. }\\
We give numerical examples by a method by first determining 
a source.  Our numerical results indicated acceptably good accuracy 
also with noisy data.   Our rather good stability can support our numerical 
results.  
\\

{\bf Acknowledgement.} 
This work was partially supported by Grant-in-Aid for Scientific Research (S) 
15H05740 of Japan Society for the Promotion of Science,
Subsidy Program for R\&D in Innovation Coast Framework,  
NSFC (No. 11771270, 91730303) and the 
\lq\lq RUDN University Program 5-100''. 
J. Cheng is supported by NSFC (No. 11421110002, No. 11771270). 
Most of the paper was completed during the stay of the fourth author at 
the Universit\`a degli Studi di Napoli Federico II in January 2019.


\begin{thebibliography}{99}%
\bibitem{CHY}
J. Cheng, Y.C. Hon and M. Yamamoto, 
Stability in line unique continuation of harmonic 
functions: general dimensions, 
J. Inv. Ill-Posed Problems {\bf 6} (1998) 319-326.

\bibitem{CPY}
J. Cheng, S. Pr\"ossdorf and M. Yamamoto,
Local estimation for an integral equation of first kind with analytic kernel,
J. Inv. Ill-Posed Problems {\bf 6} (1998) 115-126.

\bibitem{CY1}
J. Cheng and M. Yamamoto, Unique continuation on a line for harmonic 
functions. Inverse Problems  {\bf 14} (1998) 869-882. 

\bibitem{CY2}
J. Cheng and M. Yamamoto, Conditional stabilizing estimation for an 
integral equation of first kind with analytic kernel, J. Integral Equations 
Appl. {\bf 12} (2000) 39-61.

\bibitem{CY3}
J. Cheng, M. Yamamoto. One new strategy for a priori
choice of regularizing parameters in Tikhonov's regularization,
Inverse Problems {\bf 16} (2000) L31-L38.

\bibitem{I}
V. Isakov, Inverse Problems for Partial Differential Equations,
Springer-Verlag, Berlin, 2006.

\bibitem{LRS}
M. M. Lavrent'ev, V. G. Romanov and S. P. Shishatskii, 
Ill-Posed Problems of Mathematical Physics and Analysis,
American Mathematical Society, Providence, Rhode Island, 1986.

\bibitem{MOMTS}
A. Malins, M. Okumura, M. Machida, H. Takemiya and K. Saito,
Fields of view for environmental radioactivity, 
submitted for Proceedings of the 2015 
International Symposium on Radiological Issues for Fukushima's 
Revitalized Future, 
https://arxiv.org/abs/1509.09125

\bibitem{S}
K. Saito, N. Petoussi-Henss and M. Zankl,
Calculation of the effective dose from environmental gamma
ray sources and its variation, Health Phys., {\bf 74} (1998),
698-706.



\end{thebibliography}
\end{document}